\documentclass[11pt]{article}
\usepackage{amsmath}
\usepackage{amssymb}
\usepackage{amscd}
\newcommand{\R}{\mathbb{R}}
\newcommand{\inr}[1]{\bigl< #1 \bigr>}

\newcommand{\sgn}{{\rm sgn}}
\newcommand{\E}{\mathbb{E}}

\newcommand{\eps}{\varepsilon}

\newcommand{\gp}{\gamma_2(F,\psi_2)}
\newcommand{\diamp}{{\rm diam}(F,\psi_1)}
\newcommand{\pif}[1]{\pi_{#1}(f)}

\newtheorem{Theorem}{Theorem}[section]
\newtheorem{Lemma}[Theorem]{Lemma}

\newtheorem{Definition}[Theorem]{Definition}

\newtheorem{Corollary}[Theorem]{Corollary}
\newtheorem{Remark}[Theorem]{Remark}
\newtheorem{Example}[Theorem]{Example}

\newtheorem{Question}[Theorem]{Question}
\numberwithin{equation}{section} 

\def \proof {\noindent {\bf Proof.}\ \ }

\def \endproof
{{\mbox{}\nolinebreak\hfill\rule{2mm}{2mm}\par\medbreak}}
\def\IND{1}

\begin{document}
\title{On Weakly Bounded Empirical Processes}
\author{
Shahar Mendelson\footnote{Centre for Mathematics and its
Applications, Institute of Advanced Studies, The Australian
National University, Canberra, ACT 0200, Australia, and Department
of Mathematics, Technion, I.I.T, Haifa 32000, Israel.
\newline
{\sf email: shahar.mendelson@anu.edu.au}} }

\medskip
\maketitle
\begin{abstract}
Let $F$ be a class of functions on a probability space
$(\Omega,\mu)$ and let $X_1,...,X_k$ be independent random variables
distributed according to $\mu$. We establish high probability tail
estimates of the form $\sup_{f \in F} |\{i : |f(X_i)| \geq t \}$
using a natural parameter associated with $F$. We use this result to
analyze weakly bounded empirical processes indexed by $F$ and
processes of the form $Z_f=\left|k^{-1}\sum_{i=1}^k
|f|^p(X_i)-\E|f|^p\right|$ for $p>1$. We also present some geometric
applications of this approach, based on properties of the random
operator $\Gamma=k^{-1/2}\sum_{i=1}^k \inr{X_i,\cdot}e_i$, where the
$(X_i)_{i=1}^k$ are sampled according to an isotropic, log-concave
measure on $\R^n$.
\end{abstract}

\section{Introduction}
Empirical Processes theory focuses on understanding the behavior
of the supremum of the process
$$
f \to Z_f = \left|\frac{1}{k} \sum_{i=1}^k f(X_i) -\E f \right|
$$
where $F$ is a class of functions on a probability space
$(\Omega,\mu)$, $f \in F$ and $(X_i)_{i=1}^k$ are independent
random variables distributed according to $\mu$. Let $\mu_k$
denote the random empirical measure $k^{-1}\sum_{i=1}^k
\delta_{X_i}$, and for a class $F$ we denote the supremum of the
empirical process indexed by $F$ by $\|\mu_k-\mu\|_F$. Often, one
would like to bound this supremum using geometric properties of
the set $F$, but the question we tackle here is slightly
different; our aim is to bound the supremum of the empirical
process indexed by powers of the class $F$, that is, the supremum
of the process indexed by the set $F^p \equiv \{ |f|^p : f \in
F\}$ for $p>1$ using the geometry of the set $F$ rather than the
geometry of $F^p$. The difficulty arises when elements in $F$ are
not necessarily bounded functions, or in cases where the
$L_\infty$ bound is weak - while the situation is considerably
simpler in the bounded case. For example, if $F$ consists of
functions bounded by $1$ then the empirical process indexed by
$F^p$ can be bounded using a combination of symmetrization and
contraction arguments. Indeed, by the Gin\'{e}-Zinn symmetrization
method (see, for example, \cite{GZ84,VW}),
\begin{align*}
\E \|\mu_k-\mu\|_{F^p} \leq & 2\E \sup_{f \in F} \left|
\frac{1}{k} \sum_{i=1}^k \eps_i |f|^p(X_i) \right|
\\
\leq & 2p \E \sup_{f \in F} \left| \frac{1}{k} \sum_{i=1}^k \eps_i
f(X_i) \right|,
\end{align*}
where $(\eps_i)_{i=1}^k$ are independent, symmetric
$\{-1,1\}$-valued random variables. The last inequality is evident
from a contraction principle \cite{LT,VW} and the fact that $|x|^p$
is a Lipschitz function on $[-1,1]$ with constant $p$.

Moreover, for a class of uniformly bounded functions, the supremum
of the empirical process $\|\mu_k-\mu\|_F$ is highly concentrated
around its mean, as the following theorem, due to Talagrand,
shows.
\begin{Theorem} \label{thm:Talagrand} \cite{Tal94,Led}
Let $F$ be a class of mean zero functions defined on
$(\Omega,\mu)$ such that for every $f \in F$, $\|f\|_\infty \leq
b$. Let $X_1,...,X_k$ be independent random variables distributed
according to $\mu$ and set $\sigma^2=k \sup_{f \in F} {\rm
var}(f)$. Define
  \begin{equation*}
    Z  =\sup_{f \in F} \sum_{i=1}^k f(X_i), \ \ \
    \bar{Z}  = \sup_{f \in F} \left| \sum_{i=1}^k f(X_i) \right|.
  \end{equation*}
Then, for every $x>0$,
\begin{equation} \label{eq:Talagrand}
Pr \left(\left\{\left|Z-\E Z\right| \geq x \right\} \right) \leq
c_1\exp\left(-\frac{x}{c_2b}\log\left(1+\frac{bx}{\sigma^2+b\E
\bar{Z}}\right)\right),
\end{equation}
where $c_1$ and $c_2$ are absolute constants. The same inequality
is also true when $\bar{Z}$ replaces $Z$ in \eqref{eq:Talagrand}.
\end{Theorem}
Unfortunately, in many applications the function class at hand
does not consist of uniformly bounded functions, or even if the
functions are, the uniform bound is very bad. One such example
which motivated this study is the class of linear functionals of
Euclidean norm 1 on $\R^n$, and the variables $X_i$ are
distributed according to a Borel measure on $\R^n$ which is
natural from the geometric viewpoint, namely, a measure which is
isotropic and log-concave.
\begin{Definition}
A probability measure $\mu$ on $\R^n$ is called isotropic if for
every $y \in \R^n$, $\int |\inr{x,y}|^2 d\mu(x) = \|y\|^2$. The
measure $\mu$ is log-concave if for every $0<\lambda<1$ and every
Borel measurable $A,B \subset \R^n$, $\mu(\lambda A+(1-\lambda)B)
\geq \mu(A)^\lambda\mu(B)^{1-\lambda}$, where $A+B$ is the Minkowski
sum of $A$ and $B$.
\end{Definition}

A question of particular interest in this case can be formulated
as follows:
\begin{Question} \label{qu:sphere}
Let $\mu$ be an isotropic measure on $\R^n$ and let $X_1,...,X_k$
be independent, distributed according to $\mu$. Given $T \subset
\R^n$, for every $0<\eps,\delta<1$ and $p>1$, what is the smallest
integer $k_0$ such that for every $k \geq k_0$, with probability
at least $1-\delta$,
$$
\sup_{t \in T} \left|\frac{1}{k}\sum_{i=1}^k |\inr{X_i,t}|^p -\E
|\inr{X,t}|^p \right| < \eps?
$$
\end{Question}

Two simple examples which come to mind are when $p=2$, $T=S^{n-1}$
and $\mu$ is the Gaussian measure on $\R^n$ or the uniform measure
on the vertices of the unit cube.
\begin{Example} \label{ex:gauss}
{\rm For every $t \in \R^n$ define the linear functional
$f_t=\inr{t,\cdot}$ and set $F=\{f_t:t \in S^{n-1}\}$. Let $\mu_G$
be the Gaussian measure on $\R^n$ and note that for every $t \in
S^{n-1}$, $\E f_t^2 =1$. Then,
$$
\|\mu_k-\mu\|_{F^2}=\sup_{t\in S^{n-1}} \left| \frac{1}{k}
\sum_{i=1}^k \inr{t,X_i}^2 -1 \right| = \sup_{t \in S^{n-1}}
\left| \frac{1}{k} \|\Gamma t\|^2 - 1 \right|,
$$
where $\Gamma$ is a random $k \times n$ matrix with independent,
standard Gaussian random variables as entries. Hence, if
$\|\mu_k-\mu\|_{F^2}<\eps$, the gaussian matrix is an almost
isometric embedding of $\ell_2^n$ in $\ell_2^k$ which is a well
known and useful fact and occurs as long as $k \geq
c(\eps,\delta)n$ (see \cite{Pis}). Another example is when
$\mu=\mu_R$ is the uniform probability measure on $\{-1,1\}^n$.
Thus, if $\|\mu_k-\mu\|_{F^2} \leq \eps$ then a random $k \times
n$ matrix with independent, symmetric, $\{-1,1\}$-valued entries
is an almost isometric embedding of $\ell_2^n$ in $\ell_2^k$.
Unfortunately, functions in $F$ on the probability space
$(\R^n,\mu_G)$ are not bounded, while on $(\R^n, \mu_R)$ the best
uniform $L_\infty$ bound is $\sup_{t \in S^{n-1}} \|f_t\|_\infty
\leq \sqrt{n}$ which is too weak to be useful. Therefore,
symmetrization and concentration methods which are so helpful in
the bounded case can not assist in resolving Question
\ref{qu:sphere} here, as well as in other, more general examples
we will explore.

The useful property of linear functionals (with respect to both
$\mu_G$ and $\mu_R$) is that for every $f_t \in F$,
$$
Pr \left(|f_t| \geq u\right) \leq 2\exp\left(-cu^2\right)
$$
for a suitable absolute constant $c$, implying that functions in
$F$ exhibit a subgaussian behavior. Moreover, using Borell's
inequality \cite{Bor,MS}, one can show that if $\mu$ is an
arbitrary isotropic log-concave measures, linear functionals
exhibit a subexponential decay.

To formulate these decay properties in a more accurate way, we
require the definition of Orlicz norms \cite{LT,VW}}.
\end{Example}
\begin{Definition} \
For $\alpha \geq 1$ the $\psi_\alpha$ norm of a random variable $Y$
is defined by
$$
\|Y\|_{\psi_\alpha} = \inf \left\{ u>0 : \E
\exp(|Y|^\alpha/u^\alpha) \leq 2 \right\}.
$$
\end{Definition}
It is standard to verify that if $Y$ has a bounded $\psi_\alpha$
norm then $Pr \left(|Y| \geq t \right) \leq 2\exp(-c
t^\alpha/\|Y\|^\alpha_{\psi_\alpha})$ where $c$ is an absolute
constant. The reverse direction is also true, and if $Y$ has a
tail bounded by $\exp(-t^\alpha/K^\alpha)$ then
$\|Y\|_{\psi_\alpha} \leq c_1 K$.

Out main goal is to show how decay properties of individual class
members can be combined to control $\|\mu_k-\mu\|_{F^p}$.

As a starting point, let us consider the linear case where is
addition, functionals are subgaussian with respect to the
$\ell_2^n$ norm, i.e. for every $y \in \R^n$,
$\|\inr{y,X}\|_{\psi_2} \leq c\|y\|_2$. In particular, the
diameter of $F=S^{n-1}$ is bounded with respect to the $\psi_2$
norm.

This fact by itself is not enough to bound $\|\mu_k - \mu
\|_{F^p}$, and to that end we require the following notion of
complexity of the class $F$.
\begin{Definition} \label{def:gamma-2} \cite{Tal:book}
For a metric space $(T,d)$, an {\it admissible sequence} of $T$ is a
collection of subsets of $T$, $\{T_s : s \geq 0\}$, such that for
every $s \geq 1$, $|T_s|=2^{2^s}$ and $|T_0|=1$. For $\beta \geq 1$,
define the $\gamma_\beta$ functional by
$$
\gamma_\beta(T,d) =\inf \sup_{t \in T} \sum_{s=0}^\infty
2^{s/\beta}d(t,T_s),
$$
where the infimum is taken with respect to all admissible sequences
of $T$.
\end{Definition}

In \cite{MPT} the question of estimating $\|\mu_k-\mu\|_{F^2}$ has
been studied for sets of functions which have a bounded diameter
with respect to the $\psi_2$ metric and a finite
$\gamma_2(F,\psi_2)$, under the additional assumption that for
every $f \in F$, $\E f^2 =1$.
\begin{Theorem} \label{thm:MPT} \cite{MPT}
There exist absolute constants $c_1, c_2, c_3$ and for which the
following holds. Let $(\Omega,\mu)$ be a probability space, set
$F$ to be a subset of the unit sphere of $L_2(\mu)$ and assume
that ${\rm diam}(F, \psi_2) = \alpha$. Then, for any $\theta>0$
and $k\ge 1$ satisfying
$$
c_1\alpha \gamma_2(F, \psi_2)\leq \theta\sqrt{k},
$$
with probability at least $1-\exp(-c_2\theta^2k/\alpha^4)$, $
\|\mu_k-\mu\|_{F^2} \leq \theta$. Moreover, if $F$ is symmetric,
then $\E \|\mu_k-\mu\|_{F^2} \leq c_3\alpha
\gamma_2(F,\psi_2)/\sqrt{k}$.
\end{Theorem}
Theorem \ref{thm:MPT} gives an answer to Question \ref{qu:sphere}
for $p=2$ under a $\psi_2$ assumption in a very general situation.
It is particularly helpful when the $\psi_2$ metric endowed on $F$
is equivalent to the $L_2$ metric, that is, if for every $f,g \in
F$, $\|f-g\|_{\psi_2} \leq K \|f-g\|_{L_2}$. In such a case, ${\rm
diam}(F,\psi_2) \sim {\rm diam}(F,L_2)$ and $\gamma_2(F,\psi_2)
\sim \gamma_2(F,L_2)$, where by $A \sim B$ we mean that there are
absolute constants $c$ and $C$ such that $cA \leq B \leq CA$. By
the majorizing measures Theorem (see \cite{Tal:book} for the most
recent survey on the subject), $\gamma_2(F,L_2)$ is equivalent to
the expectation of the supremum of the Gaussian processes indexed
by $F$, denoted by $\E\|G\|_F$. Therefore, under a $\psi_2$
assumption, Theorem \ref{thm:MPT} implies that if $F \subset
S(L_2)$ then for every $0<\delta<1$, with probability at least
$1-\delta$,
$$
\|\mu_k-\mu\|_{F^2} \leq c \frac{\E\|G\|_F}{\sqrt{k}},
$$
where $c$ depends on $\delta$ and on the equivalence constant
between the $\psi_2$ and $L_2$ metrics.

In the geometric context of Example \ref{ex:gauss}, Theorem
\ref{thm:MPT} is helpful when the indexing set in an arbitrary
subset of $S^{n-1}$. Moreover, if the measure $\mu$ happens to be
isotropic, then the Gaussian process indexed by $F$ is the
isonormal one and thus $\gamma_2(F,L_2) \sim \E\sup_{t \in T}
\left|\sum_{i=1}^n g_it_i \right|$, where $g_1,..,g_n$ are
independent, standard Gaussian variables.


Unfortunately, the assumption that the $\psi_2$ metric is
equivalent to the $L_2$ metric is overly optimistic. In
particular, the class may not have a well bounded diameter in
$\psi_2$, or the diameter could be of the same order of magnitude
as $\gamma_2(F,\psi_2)$. For example,
if $\mu$ is log-concave and isotropic, then for every $y \in
\R^n$, the function $f_y=\inr{y,\cdot}$ satisfies
$\|f_y\|_{\psi_1(\mu)} \leq K \|f_y\|_{L_2(\mu)}$ and the $\psi_1$
and $L_2$ norms are equivalent on $\R^n$, but in contrast, the
$\psi_2$ diameter of $S^{n-1}$ might be polynomial in the
dimension (e.g. $\sqrt{n}$ when $\mu$ is the normalized volume
measure on the isotropic position of the unit ball of $\ell_1^n$).
Hence, the bound one can establish from Theorem \ref{thm:MPT} is
useless in such cases because of the way it depends on the
$\psi_2$ diameter of the set.

It would be desirable to prove a result of a similar flavor to
Theorem \ref{thm:MPT}, with the $\psi_2$ diameter of $F$ replaced
by the $\psi_1$ diameter and also removes the restrictions that
$p=2$ and that $T \subset S(L_2)$. Our main result implies just
that.

To see why the $\psi_1$ case is considerably more difficult than
the $\psi_2$ one, consider a single function $h \in L_{\psi_1}$.
By Bernstein's inequality (Lemma \ref{lemma:Bern} below),
empirical means of $h$ are highly concentrated around its
expectation, with a tail which decays exponentially in sample
size. Clearly, if a function $f \in L_{\psi_2}$ then $f^2 \in
L_{\psi_1}$ and hence exhibits the degree of concentration needed
in the proof of Theorem \ref{thm:MPT}. On the other hand, if $f
\in L_{\psi_1}$ the degree of concentration of empirical means of
$f^2$ around $\E f^2$ is not strong enough for that approach.

To overcome this obstacle, the method we suggest here is to
decompose $F$ to two subsets $F_1$ and $F_2$ which satisfy that $F
\subset F_1+F_2$.





Fix $\theta(k)>0$ and consider the sets $F_1=\{f\IND_{\{|f| \leq
\theta\}} : f \in F\}$ and $F_2=\{f\IND_{\{|f|> \theta\}} : f \in
F\}$.

Since all the functions in $F_1$ are bounded by $\theta$, the
empirical mean $\mu_k(f)$ is highly concentrated around the true
mean for any $f \in F_1$ and $\|\mu_k-\mu\|_{F_1}$ (or
$\|\mu_k-\mu\|_{F_1^p}$, using a contraction argument) is well
behaved. The key point in this approach is to control the ``large
part" of the process, namely,
$$
\sup_{f \in F} k^{-1}\sum_{i=1}^k |f|^p\IND_{\{|f| > \theta \}},
$$
and to show that the supremum is small even for a relatively low
level of truncation $\theta$. The reason this supremum is small
has nothing to do with the concentration of each individual class
member around its mean, but rather with the fact that with high
probability, all the functions in $F$ have an empirical
distribution which decays quickly. And indeed, the main Theorem we
present is an ``empirical processes" version of result due to
Bourgain on the distribution of functions in $F$ with respect to
the (random) empirical measure $\mu_k$.

\noindent{\bf Theorem A.} {\it There exist absolute constants
$c_1$, $c_2$ and $c_3$ for which the following holds. Let $F$ be a
class of mean zero functions on $(\Omega,\mu)$. For every $v_1,v_2
\geq c_1$, with probability at least $1-\exp(-c_2 \min\{v_1,v_2 \}
)$, for any $f \in F$ and $t>0$,
\begin{equation*}
\left| \left\{ i : \left|f(X_i)\right| \geq t \right\} \right|
\leq \max \left\{ \frac{c_3v_1 \gamma_2^2(F,\psi_2)}{t^2}, ek\exp
\left(-\frac{\theta}{c_3\alpha v_2}\right) \right\},
\end{equation*}
where $\alpha = \diamp$. }

Bourgain's argument \cite{Bour} is very different from ours and is
tailored to the specific case $F=\{\inr{y,\cdot} : y \in
S^{n-1}\}$, where $X_1,...,X_k$ are selected according to a
log-concave measure on $\R^n$ (see Section \ref{sec:main} for a
more detailed discussion).

The proof of Theorem A is based on the following estimate (which
will be shown to be optimal) on the $\ell_1$ structure of a random
coordinate projection of $F$.

\noindent{\bf Theorem B.} {\it For every $0<\delta<1$ there is a
constant $c(\delta)$ for which the following holds. For every
integer $k$, with probability at least $1-\delta$, for every $f \in
F$ and $I \subset \{1,...,k\}$,
\begin{align*}
\sum_{i \in I} |f(X_i)| & \leq c(\delta)
\left(\sqrt{|I|}\gamma_2(F,\psi_2) + {\rm
diam}(F,\psi_1)|I|\log\left(\frac{ek}{|I|}\right)\right),
\\
\sum_{i \in I} |f(X_i)| & \leq c(\delta)
\left(\sqrt{|I|}\gamma_2(F,\psi_2) + {\rm
diam}(F,\psi_2)|I|\sqrt{\log\left(\frac{ek}{|I|}\right)}\right).
\end{align*}
}

We present several geometric applications of Theorem A. The first
of which is a ``log-concave" version of the celebrated result of
Pajor and Tomczak-Jaegermann \cite{PT} on sections of small
diameter of a convex, symmetric body $K$ (see also
\cite{MiPi,Mi,MiPa} for results along the same lines). We show
that if $X_1,...,X_k$ are selected according to an isotropic
log-concave measure on $\R^n$, then with high probability, the
intersection of the kernel of the operator $\Gamma=\sum_{i=1}^k
\inr{X_i,\cdot}e_i$ with $K$ will have a small diameter.

\noindent {\bf Theorem C.} {\it For every $0<\delta<1$ there
exists a constant $c(\delta)$ for which the following holds. Let
$\mu$ be an isotropic, log-concave measure on $\R^n$ and let $K
\subset \R^n$ be a convex symmetric body. If $X_1,...,X_k$ are
independent, distributed according to $\mu$, then with probability
at least $1-\delta$,
$$
{\rm diam}({\rm ker} \Gamma \cap K) \leq q_k^*(K),
$$
where
$$
q_k^*(K) = \inf\left\{\rho>0 : \rho \geq
c(\delta)\frac{V_\rho\sqrt{\log{V_\rho}}}{\sqrt{k}}\right\},
$$
and $V_\rho=\gamma_2(K\cap \rho S^{n-1} ,\psi_2)$. }

\vskip 0.4cm

If $\mu$ is a subgaussian measure, Theorem C gives a weaker result
(by up to a factor of $\sqrt{\log{n}}$) and with a weaker
probability estimate than Theorem \ref{thm:MPT}. On the other
hand, it is applicable for a wider set of measures.

The downside of our approach is that it depends on the parameter
$\gamma_2(F,\psi_2)$ which is often hard to bound. However, as we
show, a completely $\psi_1$ version of Theorem B is not true and
one might have to use the additional structural assumptions on the
indexing set to improve our estimate. Luckily, in the case
$F=\{f_t : t \in S^{n-1}\}$, it is possible to bound
$\|\mu_k-\mu\|_{F^p}$ in a rather strong sense (though probably
suboptimal by a logarithmic factor) using a truncation of the
measure $\mu$. Let $\mu$ be a probability measure on $\R^n$, for
every integer $k$, let $X_1,...,X_k$ be independent, distributed
according to $\mu$ and set $H_k=\E \max_{1 \leq i \leq k}
\|X_i\|$. Observe that if $Y_i=X_i \IND_{\{\|X\| \leq
c_1(\delta)H_k\}}$, then with probability at least $1-\delta$,
$X_i=Y_i$ for $1 \leq i \leq k$. Thus, one can consider the
process $\|\nu_k-\nu\|_{F^p}$ instead of the original process
$\|\mu_k-\mu\|_{F^p}$. Moreover, one can show

\noindent{\bf Theorem D.} {\it There exist absolute constants
$c_1$, $c_2$ and $c_3$ for which the following holds. If $F=\{f_t
: t \in S^{n-1}\}$ then
$$
\gamma_2(F,\psi_2(\nu)) \leq c_1 H_k\sqrt{\log{n}}.
$$
} \vskip0.5cm

Note that if $\mu$ is an isotropic log-concave measure on $\R^n$ and
if $n \leq k \leq \exp(c_2\sqrt{n})$ then $H_k \leq c_3 \sqrt{n}$,
which is a fact recently proved by Paouris \cite{Pao}.

As we demonstrate in Section \ref{sec:applications}, the
combination of Theorem B and Theorem D allows us to bound
$$
\sup_{t \in S^{n-1}} \left| \frac{1}{k}\sum_{i=1}^k
|\inr{t,X_i}|^p - \E|\inr{t,X}|^p \right|
$$
for any log-concave measure.



\section{Preliminary Results}
In this section we present basic results which are used throughout
this article. First, a notational convention. All absolute constants
are positive numbers, denoted by $c,c_1,c_2,..$ etc. Their value may
change from line to line. We denote the Euclidean norm by $\| \ \|$,
while all other norms will be clearly specified.

There are several useful results regarding the concentration and
tail behavior of sums of independent random variables. The first
one we present here deals with subgaussian random variables and
can be easily seen using the moment generating function.

\begin{Lemma} \label{lemma:sum-subgauss} \cite{VW}
There exists an absolute constant $c$ for which the following holds.
Let $X$ be a subgaussian random variable and let $X_1,...,X_k$ be
independent, distributed as $X$. Then, for every $a=(a_1,...,a_k)
\in \R^k$
$$
\|\sum_{i=1}^k a_i X_i \|_{\psi_2} \leq c\|X\|_{\psi_2} \|a\|.
$$
\end{Lemma}

If $X$ is not a $\psi_2$ random variable and only exhibits a
subexponential tail then Bernstein's inequality describes the way
the average of independent copies of $X$ concentrate around their
mean - with a tail which is a mixture of subgaussian and
subexponential.
\begin{Lemma} \label{lemma:Bern} \cite{VW}
There exists an absolute constant $c$ for which the following holds.
Let $X_1,...,X_k$ be independent copies of a mean zero random
variable. Then, for any  $t>0$,
\begin{equation*}
\Pr \left(\left|\frac{1}{k}\sum_{i=i}^k X_i \right| >t \right) \leq
2\exp \left(-c\, k\min \left( \frac{t}{\|X\|_{\psi_1}},\,
\frac{t^2}{\|X\|^2_{\psi_1}}\right)\right).
\end{equation*}
\end{Lemma}

%

It turns out that using the generic chaining method
\cite{Tal:book} combined with Lemma  \ref{lemma:sum-subgauss} or
Lemma \ref{lemma:Bern}, one can bound the supremum of the
empirical process indexed by $F$.

\begin{Theorem} \cite{Tal:book} \label{thm:generic-chaining}
There exists an absolute constant $c$ for which the following holds.
If $F$ is a class of functions on $(\Omega,\mu)$, then for every
integer $k$,
\begin{align*}
\E \|\mu_k-\mu\|_F & \leq c \frac{\gamma_2(F,\psi_2)}{\sqrt{k}},
\\
\E \|\mu_k-\mu\|_F & \leq
c\left(\frac{\gamma_2(F,\psi_1)}{\sqrt{k}}+
\frac{\gamma_1(F,\psi_1)}{k} \right),
\end{align*}
and similar bounds hold with high probability.
\end{Theorem}

In many cases, computing the $\gamma$ functionals is a difficult
task. It is possible to upper bound them using a metric entropy
integral, similar to Dudley's integral in the context of Gaussian
process.
\begin{Definition}
Let $(T,d)$ be a metric space. The covering number of $T$ at scale
$\eps$ is the minimal number of open balls (with respect to the
metric $d$) of radius $\eps$ needed to cover $T$. The covering
numbers of $(T,d)$ are denoted by $N(\eps,T,d)$.
\end{Definition}

Since one way of forming an admissible sequence for $(T,d)$ is to
use an almost optimal cover (the set $T_s$ is a cover at the scale
at which one needs $2^{2^s}$ balls to cover $T$), the following is
evident:
\begin{Lemma}
There exists an absolute constant $c$ for which the following holds.
Let $(T,d)$ be a metric space. Then,
$$
\gamma_2(T,d) \leq c\int_0^\infty \sqrt{\log{N(\eps,T,d)}}d\eps.
$$
\end{Lemma}
A much more difficult result, due to Talagrand
\cite{Tal94,Tal:book}, is that if $T$ is a unit ball of a $2$-convex
normed space, $\gamma_2$ could be bounded from above by a sharper
version of the entropy integral.

\begin{Definition}
A Banach space is called 2-convex if there is $\rho>0$ such that for
$\|x\|,\|y\| \leq 1$, $\|x+y\| \leq 2-2\rho\|x-y\|^2$.
\end{Definition}

\begin{Theorem} \label{thm:talagrand} \cite{Tal1}
For every $\rho>0$ there exists a constant $c(\rho)$ for which the
following holds. If $Y$ is a $2$-convex Banach space with
parameter $\rho$ and if the metric $d$ is given by some other norm
$| \ |$, then
$$
\gamma_2(B_Y, d) \leq c(\rho) \left(\int_0^\infty \eps \log N
\left(B_Y, B_{| \ |},\eps \right) d\eps \right)^{\frac{1}{2}}.
$$
\end{Theorem}
Theorem \ref{thm:talagrand} is used in the case $Y=\ell_2^n$, the
$n$-dimensional Euclidean space, where $d$ is the metric endowed
on $\R^n$ by the $\psi_2$ norm (see Section
\ref{sec:applications}).

\section{Decomposing classes of functions}
\label{sec:main}
Let $F$ be a class of functions on the
probability space $(\Omega,\mu)$ and assume that for every $f \in
F$, $\E f =0$.

Let us formulate the main technical tool we require.
\begin{Theorem} \label{thm:sets}
There exists absolute constants $c_1$ and $c_2$ for which the
following holds. Let $F$ be a class of mean zero functions on
$(\Omega,\mu)$ and set $X_1,...,X_k$ to be independent random
variables distributed according to $\mu$. Then, for every $v_1, v_2
\geq c_1$, with probability at least $1-\exp(-c_2
\min\{v_1^2,v_2\})$, for every $I \subset \{1,...,k\}$,
$$
\sup_{f \in F} \left| \sum_{i \in I} f(X_i) \right| \leq v_1
\sqrt{|I|} \gp + v_2 \diamp |I| \log\left(\frac{ek}{|I|}\right).
$$
\end{Theorem}

Theorem \ref{thm:sets} has a similar version in which one assumes
that the set of functions is well bounded in $\psi_2$.
\begin{Theorem} \label{thm:sets-psi_2}
There exist absolute constants $c_1$ and $c_2$ for which the
following holds. Let $F$ and $X_1,...,X_k$ be as in Theorem
\ref{thm:sets}. Then, for every $v \geq c_1$, with probability at
least $1-\exp(-c_2 v^2)$, for every $I \subset \{1,...,k\}$,
$$
\sup_{f \in F} \left| \sum_{i \in I} f(X_i) \right| \leq v
\left(\sqrt{|I|} \gp + {\rm diam}(F,\psi_2)|I|  \sqrt{
\log\left(\frac{ek}{|I|}\right)}\right).
$$
\end{Theorem}

Theorem \ref{thm:sets} is an empirical processes version of a
lemma due to Bourgain (\cite{Bour}, see also \cite{GiaMil}) which
deals with the case when $F$ is $S^{n-1}$, considered as a class
of linear functionals on $\R^n$ and $\mu$ is an isotropic
log-concave measure. Unlike Bourgain's argument, which relies
heavily on the fact that the functions in the class are linear
functionals and on that the indexing set is the whole sphere,
Theorem \ref{thm:sets} is very general.

Observe that if the $L_2$ and $\psi_2$ metrics are equivalent on
$F$ with a constant $\beta$ and if $\E \|G\|_F$ denotes the
expectation of the supremum of the Gaussian process indexed by
$F$, then by the majorizing measures Theorem there are absolute
constants $c$ and $C$ and a constant $c_1(\beta)$ depending only
on $\beta$ such that
$$
c_1(\beta) \gamma_2(F,\psi_2) \leq c\gamma_2(F,L_2) \leq \E\|G\|_F
\leq C\gamma_2(F,L_2) \leq C \gamma_2(F,\psi_2).
$$
Therefore, by Theorem \ref{thm:sets}, with probability at least
$1-\delta$, for every $I \subset \{1,...,k\}$,
\begin{align*}
\sup_{f \in F} \left| \sum_{i \in I} f(X_i) \right| & \leq
c(\delta,\beta) \left( \sqrt{|I|} \E\|G\|_F + \diamp |I|
\log\left(\frac{ek}{|I|}\right)\right),
\\
\sup_{f \in F} \left| \sum_{i \in I} f(X_i) \right| & \leq
c(\delta,\beta) \left(\sqrt{|I|} \E\|G\|_F + {\rm
diam}(F,\psi_2)|I| \sqrt{\log\left(\frac{ek}{|I|}\right)}\right).
\end{align*}

Let us point out that it is impossible to obtain a fully $\psi_1$
version of Theorem \ref{thm:sets}. Indeed, suppose the converse
was true, and that for every set $F$ and integer $k$, with
probability at least $1-\delta$, for every $I \subset
\{1,...,k\}$,

\begin{equation} \label{eq:wrong}
\sup_{f \in F} \left| \sum_{i \in I} f(X_i) \right| \leq c(\delta)
\left(\sqrt{|I|} \gamma_2(F,\psi_1) + \diamp |I|
\log\left(\frac{ek}{|I|}\right)\right).
\end{equation}
Let $Y$ be an exponential random variable and let
$X^{(n)}=\sum_{i=1}^n \frac{Y_i}{\sqrt{\log{(i+1)}}}e_i \in \R^n$
where $(e_i)_{i=1}^n$ is the standard basis in $\R^n$ and
$(Y_i)_{i=1}^n$ are independent copies of $Y$. Setting $\mu^{(n)}$
to be the measure on $\R^n$ which endows $X^{(n)}$,
\eqref{eq:wrong} can not be true for $\mu^{(n)}$ and $F_n=B_1^n$,
the unit ball in $\ell_1^n$ even when $k=1$. Indeed, using
Borell's inequality (see, e.g. \cite{MS}) or by a direct
computation as in \cite{BGMN}, it is evident the for every $t \in
\R^n$,
$$\left\|\inr{t,X^{(n)}}\right\|_{\psi_1}
\leq c\left\|\inr{t,X^{(n)}}\right\|_{L_2}=c\left(\sum_{i=1}^n
\frac{t_i^2}{\log{(i+1)}}\right)^{1/2} \equiv c|t|^{(n)},
$$
where $| \ |^{(n)}$ is the weighted Euclidean norm with weights
$(\sqrt{\log(i+1)})_{i=1}^n$ and $c$ is an absolute constant.
Hence, by the majorizing measures Theorem and a standard
computation, there are absolute constants $c$, $c_1$ and $c_2$
such that for every $n$,
$$
\gamma_2\left(F_n,\psi_1(\mu^{(n)})\right) \leq c\gamma_2(F_n,| \
|^{(n)}) \leq c_1\E \sup_{t \in B_1^n} \left|\sum_{i=1}^n g_i
\frac{t_i}{\sqrt{\log(i+1)}} \right| \leq c_2.
$$
Therefore, if \eqref{eq:wrong} were correct for $k=1$, it would
follow that with probability of at least $1/2$, $\sup_{t \in
B_1^n} \inr{t,X^{(n)}} \leq c_3$, for a suitable $c_3$ which is
independent of $n$. On the other hand, an easy computation shows
that with probability larger than some constant $c_4$,
$$
\sup_{t \in B_1^n} \inr{t,X^{(n)}} \geq \sqrt{\log{(n+1)}},
$$
and thus it is impossible to get a completely $\psi_1$ version of
Theorem \ref{thm:sets}.

Next, observe that Theorem \ref{thm:sets-psi_2} is optimal, in the
sense that both the $\gamma_2$ term and the term that depends on
the $\psi_2$-diameter are required. To see this, fix an integer
$k$ and let $1 \leq \ell <k$. Set $F=\{a,-a\} \subset S^{n-1}$,
acting as linear functional of $\R^n$ and let $X=(g_1,...,g_n)$ be
a Gaussian vector in $\R^n$. With this choice of $X$, the $\psi_2$
and $\ell_2$ metrics on $\R^n$ are equivalent with an absolute
constant, and thus $\gamma_2(F,\psi_2) \leq c \gamma_2(F,\ell_2)
=c_1$. On the other hand, writing $a=(a_1,...,a_n)$, for every $1
\leq \ell \leq k$,
\begin{equation*}
\sup_{f \in F} \sup_{I \in E_\ell} \left|\sum_{i \in I}
f(X_i)\right|=\sup_{I \in E_\ell} \left|\sum_{i \in I} \sum_{j=1}^n
a_jg_{i,j} \right|,
\end{equation*}
which is the supremum of the Gaussian process indexed by
$$
E_\ell=\left\{I\ : \ I \subset \left\{1,...,k\right\}, \ |I|=\ell
\right\}
$$
with the covariance structure endowed by the Hamming metric on
$E_\ell$, given by $d_H(I,J)=|I \vartriangle J|^{1/2}$. Recall the
well known entropy estimate for $E_\ell$ with respect to this
metric (see, for example, \cite{MPR}):

\begin{Lemma} \label{Lemma:MPR}
For $0<\lambda\le 1/2$ there exists a constant $c_\lambda$ for which
the following holds. For every integers $k$ and $1 \leq \ell \leq
k$, there is a subset $P \subset E_\ell$ which satisfies that $\log
|P| \geq (1-\lambda) \ell \log\left(c_\lambda \frac{k}{\ell}\right)$
and if $I, J \in P$ and $I \not=J$ then $d_H(I,J) \geq \sqrt{\lambda
\ell}$. In other words,
$$
\log N \left(E_\ell,\sqrt{\lambda \ell}, d_H\right) \geq (1-\lambda)
\ell \log \left(c_\lambda \frac{k}{\ell}\right).
$$
\end{Lemma}
Combining Lemma \ref{Lemma:MPR} for $\lambda =1/4$ with Sudakov's
minoration (see, e.g. \cite{LT}), it is evident that
$$
\E \sup_{I \in E_\ell} \left|\sum_{i \in I} \sum_{j=1}^n a_jg_{i,j}
\right| \geq c\ell \sqrt{\log\left(\frac{ck}{\ell}\right)},
$$
proving that the second term in Theorem \ref{thm:sets-psi_2} is
indeed necessary.

To show that the $\gamma_2$ term is necessary, let $F=\{-1,1\}^n$
acting as linear functionals, and again set $X$ to be the Gaussian
vector on $\R^n$. Then, for every $1 \leq \ell \leq k$,
\begin{equation*}
\sup_{a \in \{-1,1\}^n} \sup_{I \in E_\ell} \left|\sum_{i \in I}
\sum_{j=1}^n g_{i,j}a_j \right| \geq \sup_{a \in \{-1,1\}^n}
\left|\sum_{i = 1}^\ell \sum_{j=1}^n g_{i,j}a_j \right|.
\end{equation*}
The latter is the supremum of the Gaussian process indexed by
$\{-1,1\}^n$ with the covariance structure given by the metric
$d(u,v)=\sqrt{\ell}\|u-v\|_{\ell_2^n}$. Thus, it is standard to
verify that
$$
\E\sup_{a \in \{-1,1\}^n} \left|\sum_{i \in 1}^\ell \sum_{j=1}^n
g_{i,j}a_j \right| \geq c\sqrt{\ell}n.
$$
On the other hand, ${\rm diam}(\{-1,1\}^n,\psi_2) \leq c{\rm
diam}(\{-1,1\}^n,\ell_2^n) \leq c\sqrt{n}$. Thus, the upper bound
from Theorem \ref{thm:sets-psi_2} is of the order of
$\sqrt{\ell}n+\sqrt{n}\sqrt{\ell \log(ek/\ell)}$, showing that the
$\gamma_2$ term can not be removed from the bound.

\noindent {\bf proof of Theorem \ref{thm:sets}.} To control $\sup_{f
\in F} \left|\sum_{i \in I} f(X_i) \right|$, consider the following
$k$ processes. Recall that for every $1 \leq \ell \leq k$, $E_\ell =
\{ I : I \subset \{1,...,k\}, \ |I| =\ell \}$ and define the random
process
$$
f \to Z_f^\ell = \sup_{I \in E_\ell} \left| \sum_{i \in I} f(X_i)
\right|,
$$
where $X_1,...,X_k$ are independent random variables distributed
according to $\mu$.

Fix $1 \leq \ell <k$ (the result for $\ell=k$ requires minor changes
and is omitted) and consider the process $Z_f^\ell$. Observe that
for every $f,g \in F$,
\begin{align*}
Pr \left( \left| Z_f^\ell - Z_g^\ell \right| \geq t \right) & \leq
Pr \left( \sup_{I \in E_\ell} \left| \sum_{i \in I} (f-g)(X_i)
\right| \geq t \right) \\
& \leq |E_\ell| Pr \left( \left| \sum_{i=1}^\ell (f-g)(X_i) \right|
\geq t \right)
\\
& \leq 2|E_\ell| \exp \left(-\frac{ct^2}{\|f-g\|^2_{\psi_2} \ell}
\right),
\end{align*}
where $c$ is an absolute constant.

Without loss of generality, assume that $\gamma_2(F,\psi_2) <
\infty$, let $(F_s)_{s \geq 0}$ be an almost optimal admissible
sequence for the metric space $(F,\psi_2)$ and set $\pi_s(f)$ to
be a nearest element to $f$ in $F_s$ with respect to the $\psi_2$
metric. Thus, $|F_s| \leq 2^{2^s}$. Define $s_0$ as the first
index such that $2^{s_0-1} < \log |E_\ell| \leq 2^{s_0}$, and note
that
$$
Z_f^\ell = Z^\ell_{\pif{s_0}} + \sum_{i=s_0}^\infty
Z^\ell_{\pif{i+1}} - Z^\ell_{\pif{i}}.
$$
Fix $u>0$ to be specified later and $s \geq s_0$, and consider $t_s
= u\sqrt{\ell} \| \pif{s+1}-\pif{s}\|_{\psi_2} 2^{s/2} \sqrt{\log
|E_\ell|}$. Then,
\begin{align*}
Pr \left( \left|Z^\ell_{\pif{s+1}} - Z^\ell_{\pif{s}} \right| \geq
t_s \right) & \leq 2 |E_\ell| \exp(-cu^2 2^{s-1} \log |E_\ell|) \\
& \leq 2 \exp \left(-c\log |E_\ell| \left(u^22^{s-1} -1 \right)
\right).
\end{align*}
Take $u=v_1 /\sqrt{\log |E_\ell|}$ for $v_1 \geq c_1$ and note
that $2^s > \log |E_\ell|$, implying that the tail is upper
bounded by $2\exp(-c_2v_1^2 2^s)$. Summing over $s_0 \leq s <
\infty$ it is evident that with probability at least
$$
1-2\sum_{s_0}^\infty \exp(-c_2v_1^2 2^s) \geq 1 -
2\exp(-c_32^{s_0}v_1^2) \geq 1 -2\exp(-c_3v_1^2 \log |E_\ell|),
$$
for every $f \in F$
\begin{align*}
\sum_{i=s_0+1}^\infty \left| Z^\ell_{\pif{s+1}} - Z^\ell_{\pif{s}}
\right| & \leq v_1 \sqrt{\ell} \sum_{i=s_0+1}^\infty 2^{s/2} \|
\pif{s+1} -\pif{s} \|_{\psi_2} \\
& \leq c_4 v_1 \sqrt{\ell} \gp.
\end{align*}
To handle $F_{s_0}=\{\pif{s_0} : f \in F \}$, note that the
cardinality of this set is at most $2^{2^{s_0}} \leq 2^{2 \log
|E_\ell|}$. Applying Bernstein's inequality (Lemma
\ref{lemma:Bern}), for every $t>0$ and every $f \in F$,
\begin{align*}
Pr \left( \left| Z_f^\ell \right| \geq t\ell \right) & \leq |E_\ell|
Pr \left( \left| \sum_{i=1}^\ell f(X_i) \right| \geq t \ell \right)
\\
& \leq 2 |E_\ell| \exp \left(-c \ell \min \left\{
\frac{t^2}{\|f\|_{\psi_1}^2}, \frac{t}{\|f\|_{\psi_1}} \right\}
\right).
\end{align*}
Let $t=\frac{\log |E_\ell|}{\ell} \|f\|_{\psi_1} v_2$ for $v_2
\geq 1$. Since $1 \leq \ell < k$ then $\ell^{-1} \log |E_\ell|
\geq 1$ and $t \geq \|f\|_{\psi_1}$. Therefore, with probability
at least
$$
1-2|E_\ell|^2 \exp \left(-c_5v_2 \log |E_\ell|\right) \geq
1-2\exp(-\log |E_\ell| (c_5v_2 - c_6) ),
$$
for every $ f \in F_{s_0}$,
$$
Z_f^\ell \leq v_2 \|f\|_{\psi_1} \log |E_\ell| \leq v_2 \diamp
\log |E_\ell|.
$$
To conclude, there are absolute constants $c_7$, $c_8$ and $c_9$
such that for every $1 \leq \ell \leq k$, if $v_1,v_2 \geq c_7$,
with probability at least $1-2\exp(-c_8 \log |E_\ell|
\min\{v_1^2,v_2\})$,
$$
\sup_{f \in F} |Z_f^\ell| \leq c_9 \left(v_1 \sqrt{\ell} \gp + v_2
\diamp \log |E_\ell|\right).
$$
Summing the probabilities, the latter holds for every $1 \leq \ell
\leq k$ with probability at least $1-\exp(-c_{10}
\min\{v_1^2,v_2\})$, completing the proof.
\endproof

The proof of Theorem \ref{thm:sets-psi_2} is similar and is omitted.

\vskip 0.2cm

\noindent{\bf Proof of Theorem B.} In the $\psi_1$ case, take
$v_1=\sqrt{\log(1/\delta)}$ and $v_2=\log(1/\delta)$ for $\delta$
small enough. Fix any $f \in F$ and for $I \in E_\ell$ let
$I^+(f)=\{i : f(X_i)>0\} \cap I$ and $I^-(f)=\{i : f(X_i)<0\} \cap
I$. Then, by Theorem \ref{thm:sets}, with probability at least
$1-\delta$,
\begin{align*}
\sum_{i \in I} |f(X_i)| & = \big| \sum_{i \in I^+(f)} f(X_i) \big|
+ \big| \sum_{i \in I^-(f)} f(X_i)\big|
\\
& \leq c(\delta) \left(\sqrt{\ell}\gamma_2(F,\psi_2)+\diamp \ell
\log\left(\frac{ek}{\ell}\right)\right),
\end{align*}
as claimed. The $\psi_2$ case is equally easy.
\endproof

For Theorem \ref{thm:sets} one can derive the following uniform
empirical tail estimate for functions in $F$, which was formulated
as Theorem A in the introduction.

\begin{Corollary} \label{cor:tail-of-sets}
There exist absolute constants $c_1$, $c_2$ and $c_3$ for which
the following holds. Let $F$ be as in Theorem \ref{thm:sets}. For
every $v_1,v_2 \geq c_1$, with probability at least $1-\exp(-c_2
\min\{v_1^2,v_2 \} )$, for any $f \in F$ and $t >0$,
\begin{equation} \label{eq:theta-est}
\left| \left\{ i : \left|f(X_i)\right| \geq t \right\} \right|
\leq \max \left\{ \frac{c_3v_1^2 \gamma_2^2(F,\psi_2)}{t^2},
ek\exp \left(-\frac{t}{c_3\alpha v_2}\right) \right\},
\end{equation}
where $\alpha = \diamp$.
\end{Corollary}

\proof Fix $v_1,v_2$ as in Theorem \ref{thm:sets} and consider the
set for which the assertion of Theorem \ref{thm:sets} holds. Let
$(X_1,...,X_k)$ be in that set and for $f \in F$ and $t>0$ put
$$
I_t(f) =\{i : |f(X_i)| \geq t \}.
$$
Setting $\alpha= \diamp$ there are two possibilities. First, if
$$
2\alpha v_2 | I_t(f)| \log \left(\frac{ek}{|I_t(f)|}\right) \leq
\frac{t}{2} |I_t(f)|,
$$
then by Theorem \ref{thm:sets},
\begin{align*}
t |I_t(f)| &  \leq 2v_1 \sqrt{|I_t(f)|} \gp + 2v_2 \alpha |I_t(f)|
\log \left(\frac{ek}{|I_t(f)|}\right)
\\
& \leq 2v_1\sqrt{|I_t(f)|} \gp  + \frac{t}{2} |I_t(f)|.
\end{align*}
Thus, $t |I_t(f)|/2 \leq 2v_1 \sqrt{|I_t(f)|} \gp$, implying that
$$
|I_t(f)| \leq 16v_1^2 \frac{\gamma_2^2(F,\psi_2)}{t^2}.
$$
Otherwise, $2\alpha v_2 |I_t(f)| \log(ek/|I_t(f)|) \geq t
|I_t(f)|/2$, or in other words,
$$
|I_t(f)| \leq ek \exp\left(-\frac{t}{4v_2 \diamp}\right).
$$
\endproof
Now we are ready to formulate and prove the main theorem of this
section, which is a decomposition result for the class $F$.
\begin{Theorem} \label{thm:main}
There exist absolute constants $c_1$, $c_2$ and $c_3$, and for every
$1 \leq p<\infty$ there exists a constants $c_4(p)$ for which the
following holds. Let $F$ be a class of mean zero functions. For $v
\geq c_1$, $A \geq \gamma_2(F,\psi_2)$, $B \geq {\rm
diam}(F,\psi_1)$ and an integer $k$, set
$$
\theta \geq \max\left\{c_2vB \log
\left(c_2\frac{B^2k}{A^2}v+1\right), c_2pB\log(c_2pB+1) \right\}.
$$
Then, there are Lipschitz functions $\phi:\R \to \R$ and $\psi:\R
\to \R$ which depend on $\theta$, such that $\|\phi\|_{{\rm lip}},
\ \|\psi\|_{{\rm lip}} \leq 1$ and setting $F_{1}=\{\phi(f) : f
\in F\}$ and $F_{2}=\{\psi(f) : f \in F\}$,
\begin{description}
\item{1.} $F \subset F_{1} + F_{2}$.
\item{2.} For every $h \in F_{1}$, $\|h\|_{\infty} \leq \theta$
and for every $h \in F_2$, $\E|h|^p \leq A^2/k$ .
\item{3.} With probability at least $1-\exp(-c_3v)$,
$$
\sup_{h \in F_{2}} \sum_{i=1}^k |h(X_i)|^p \leq
c_2vA^2\left(\theta^{p-2}+\kappa_p\right),
$$
where $\kappa_p =c_4(p)\theta^{p-2}$ for $p<2$,
$\kappa_2=c_4(2)\log A$, while for $p>2$,
$\kappa_p=c_4(p)A^{p-2}$.
\end{description}
\end{Theorem}

Theorem \ref{thm:main} implies that $F$ can be decomposed into two
simple sets $F_{1}$ and $F_{2}$ (which depend on $k,p$ and $v$).
The fact that these sets are as simple as $F$ is evident because
they are images of $F$ via Lipschitz functions with constant $1$.
In particular, $\gamma_\beta(F_{i},d) \leq \gamma_\beta(F,d)$ with
respect to any reasonable metric $d$. The sets $F_{1}$ and $F_{2}$
have additional properties. $F_{1}$ has a bounded diameter in
$L_\infty$ - up to a logarithmic term, its diameter in $L_\infty$
is proportional to the $\psi_1$ diameter of $F$. Thus, if $F$ has
a well bounded diameter with respect to the $\psi_1$ metric then
functions in $F_{1}$ are highly concentrated around their means,
and one can safely use a contraction argument when bounding the
empirical process indexed by a power of $F_{1}$. The main
difficulty is in controlling the ``large part" of $F$ - i.e.
$F_{2}$. The empirical process indexed by $F_{2}^p$ is small not
because of concentration, but because the $\ell_p^k$ diameter of a
random coordinate projection of $F_{2}$ and its $L_p$ diameter are
small.

\vskip0.2cm

\noindent {\bf Proof of Theorem \ref{thm:main}.} Fix an integer
$k$ and $v$ for which Corollary \ref{cor:tail-of-sets} holds. The
first step is to select the Lipschitz functions $\phi$ and $\psi$;
those are simply truncation functions at the level $\theta$. For
$f \in F$, set $\phi(f)=\sgn(f)\min\{|f|,\theta\}$ and
$\psi(f)=f-\phi(f)$. Clearly, both functions have Lipschitz
constant $1$, $F \subset \phi(F)+\psi(F)$, and for $p>1$, because
$\phi(f)$ and $\psi(f)$ are supported on disjoint sets,
\begin{align*}
|f|^p & = \min\{|f|^p,\theta^p\}+\left(|f|^p - \theta^p\right)
\IND_{\left\{|f| \geq \theta \right\}}
\\
& \leq |\phi(f)|^p + |f|^p\IND_{\{|f| \geq \theta\}}.
\end{align*}

Let $A \geq \gamma_2(F,\psi_2)$ and $B \geq {\rm diam}(F,\psi_1)$.
It is evident that
\begin{equation} \label{eq:p-moment}
\E |f|^p \IND_{\{|f| \geq \theta\}} \leq  c_1(2pB)^p
\exp\left(-\frac{\theta}{c_2B}\right) \leq \frac{A^2}{k},
\end{equation}
as long as
\begin{equation} \label{eq:theta0a}
\theta \geq (c_3pB)\log (c_3pB) + c_3B \log
\left(\frac{B^2k}{A^2}\right),
\end{equation}
which is satisfied by our choice of $\theta$. Thus (2) is
established.

Turning to (3), recall that for every $t>0$ and $f \in F$,
$I_t(f)=\{i: |f(X_i)| \geq t \}$. By our choice of $v$, with
probability at least $1-\exp(-c_4v)$, for every $f \in F$, for
every $t>0$
\begin{equation} \label{eq:I(f)}
|I_t(f)| \leq \max \left\{ \frac{c_5v A^2}{t^2}, k\exp
\left(-\frac{t}{c_5 B v}\right) \right\}.
\end{equation}
Therefore, if
\begin{equation} \label{eq:theta0b}
t \geq t_0 = c_6 B v \log(c_6B v \sqrt{k}/A)
\end{equation}
for a suitable absolute constant $c_6$, the first term in
\eqref{eq:I(f)} is dominant. Note that if $t \geq
\max\{t_0,\sqrt{c_5v} A\}$ then $|I_t(f)|=0$, and since $\theta
\geq t_0$ then by a standard integration argument with respect to
the random empirical measure $\mu_k$, with probability at least
$1-\exp(-c_4v)$, for every $f \in F$
\begin{align*}
\E_{\mu_k} |f|^p \IND_{\{|f| \geq \theta\}} & \leq \theta^p
Pr_{\mu_k} \left(|f| \geq \theta \right) +
\int_{\theta}^{\sqrt{c_5v} A} pt^{p-1} Pr_{\mu_k} \left(|f| >
t\right) dt
\\
& \leq c_6 \frac{v A^2}{k}\left( \theta^{p-2} +
\int_{\theta}^{\sqrt{c_5v}A} pt^{p-3} dt \right),
\end{align*}
for which the claim follows.
\endproof

\begin{Remark} \label{rem:F_{1,k}}
Note that Theorem \ref{thm:main} enables one to bound
$\|\mu_k-\mu\|_{F_{1}}$ and thus $\|\mu_k-\mu\|_{F}$. Indeed,
pointwise, $\left|\left(\phi(f)\right)^p  - \left(\phi(g)\right)^p
\right| \leq 2 p \theta^{p-1} |f-g|$, implying that
$$
\gamma_2\left(\left(\phi(F)\right)^p,\psi_2\right) \leq cp
\theta^{p-1} \gp.
$$
By a standard generic chaining argument (see Theorem
\ref{thm:generic-chaining} and \cite{Tal:book}), for every $v>0$,
with probability at least $1-\exp(-c_1v^2)$,
\begin{align*}
\sup_{f \in F} \left| \frac{1}{k} \sum_{i=1}^k
\left(\phi(f)\right)^p(X_i) -\E \left(\phi(f)\right)^p \right| \leq
& c_2 v \frac{\gamma_2((\phi(F))^p,\psi_2)}{\sqrt{k}}
\\
\leq & c_3p\theta^{p-1}v \frac{\gp}{\sqrt{k}}.
\end{align*}
\end{Remark}

Combining this with Theorem \ref{thm:main}, it follows that with
probability at least $1-2\exp(-c_1v)$,
$$
\|\mu_k-\mu\|_F \leq c_2 v\left(p \theta^{p-1}
\frac{\gamma_2(F,\psi_2)}{\sqrt{k}}+
\frac{\gamma_2^2(F,\psi_2)}{k}
\left(\theta^{p-2}+\kappa_p+1\right)\right).
$$

\begin{Remark} \label{rem:p-momnet} Observe that by \eqref{eq:p-moment}, $\sup_{h \in
F_2} \E |h|^p \leq (cpB)^p \exp\left(-\frac{\theta}{cB}\right)$, a
fact we shall use below.
\end{Remark}

We end this section with another observation which follows easily
from the proof of Theorem \ref{thm:main}. To avoid complications, we
will formulate it only is the case we need it, which is when $F$ is
a class of linear functionals on $\R^n$ and $\mu$ is a measure on
$\R^n$. Consider the random variable $U=\sup_{f \in F} |\inr{f,X}|$,
and for every integer $k$ set $H_k = \E \max_{1 \leq i \leq k} U_i$,
where $(U_i)_{i=1}^k$ are independent copies of $U$.

\begin{Theorem} \label{thm:main-trunk}
For every $p \geq 1$ and $0<\delta,\eps<1$ there are constants
$c_1(\delta,\eps,p)$, $c_2(\delta,p)$ and $c_3(p)$ for which the
following holds. Let $F$ and $\mu$ be as above, consider the
random variable $Y^k=X\IND_{\{U \leq c_1(\delta,\eps,p)H_k\}}$ and
let $\nu$ be the probability measure on $\R^n$ corresponding to
$Y^k$. If $A \geq \gamma_2(F,\psi_2(\nu))$ and $B \geq {\rm
diam}(F,\psi_1(\nu))$, then with probability at least $1-\delta$,
\begin{equation} \label{eq:main-trunk}
\|\mu_k-\mu\|_{F^p} \leq c_2
\left(\theta^{p-1}\frac{\gamma_2(F,\nu)}{\sqrt{k}}+\frac{\gamma_2^2(F,\nu)}{k}\left(\theta^{p-2}
+ \tilde{\kappa}_p \right) \right)+c_3B^{1/2}\eps,
\end{equation}
where
$$
\theta=\max\left\{c_2B \log\left(c_2\frac{k
B^2}{A^2}+1\right),c_2pB\log(c_2pB+1)\right\},
$$
and $\tilde{\kappa}_p=1$ for $1 \leq p < 2$,
$\tilde{\kappa}_2=\log{H_k}$ and $\tilde{\kappa}_p=H_k^{p-2}$ for
$p>2$.
\end{Theorem}
Because the proof is based on the same arguments used in Theorem
\ref{thm:main} and Remark \ref{rem:F_{1,k}}, we will only give a
brief sketch of the required modifications which are that with
high probability, $X_i=Y_i$ for $1 \leq i \leq k$ and that by the
Cauchy-Schwarz inequality,
$$
\sup_{f \in F} \left| \E_\mu |f|^p -\E_{\nu}|f|^p \right| =\sup_{f
\in F} \E |f|^p \IND_{\{U \geq c_1(\delta,\eps,p)H_k\}} \leq
c_3(p)B^{1/2}\eps,
$$
for the right choice of constants. Thus, one can replace the
measure $\mu$ with the measure $\nu$ and consider the empirical
process $\|\nu-\nu_k\|_{F^p}$ instead of $\|\mu-\mu_k\|_{F^p}$.

The advantage of using the measure $\nu$ is that it is a truncated
version of $\mu$ at the ``correct" level for $F$ and the sample
size $k$. This truncation enables us to bound
$\gamma_2(S^{n-1},\nu)$, where $\nu$ is a truncation of an
isotropic, log-concave measure on $\R^n$.

\section{Applications} \label{sec:applications}

The first geometric application we present deals with sections of
small diameter of a convex, symmetric body. Let
$\Gamma=\sum_{i=1}^k \inr{X_i,\cdot} e_i$, where $X_1,...,X_k$ are
selected according an isotropic log concave measure. As we show
below, if $K$ is a convex, symmetric body in $\R^n$, then with
high probability, the diameter of $K \cap {\rm ker}(\Gamma)$ is
small. This extends a celebrated result of Pajor and
Tomczak-Jaegermann \cite{PT} which was proved in the case where
the random subspace was selected according to the Haar measure on
the Grassmann manifold ${\cal G}(n,k)$, but the same proof works
in the Gaussian case. Various versions and extensions of this
result may be found, for example, \cite{MiPi,Mi,MiPa}.

The following theorem is a formulation of version of this result
for a general $\psi_2$ operator (see \cite{MPT}). Let us introduce
the following notation: for a set $T \subset \R^n$ we denote by
$\ell_*(T)=\E \sup_{t \in T} \left|\sum_{i=1}^n g_i t_i \right|$,
the expectation of the supremum of the gaussian process indexed by
$T$. Recall that by the majorizing measures Theorem
\cite{Tal:book}, there are absolute constants $c_1$ and $c_2$ such
that for every $T \subset \R^n$,
\begin{equation} \label{eq:MM}
c_1 \gamma_2(T,\| \ \|) \leq \ell_*(T) \leq c_2 \gamma_2(T,\| \ \|).
\end{equation}
\begin{Theorem} \label{thm:MPT1} \cite{MPT}
There exists a absolute constant $c$ and $c_1$ for which the
following holds. Let $X_1,...,X_k$ be distributed according to an
isotropic measure $\mu$ on $\R^n$ and assume that for every $t \in
\R^n$, $\|\inr{t,\cdot}\|_{\psi_2} \leq \alpha \|t\|$ for some
$\alpha \geq 1$. If $K \subset \R^n$ is a convex symmetric body
then with probability at least $1-\exp(-c_1k/\alpha^4)$,
$$
{\rm diam}({\rm ker} \Gamma \cap K) \leq r_k^*(K),
$$
where
$$
r_k^*(K) = \inf\left\{\rho>0 : \rho \geq c_1\alpha^2\ell_*(K\cap
\rho S^{n-1})/\sqrt{k} \right\}.
$$
\end{Theorem}

Our result is similar (though with a weaker estimate) to Theorem
\ref{thm:MPT1}. Other than the different ways of estimating the
empirical process $\|\mu_k-\mu\|_{F^2}$, the two proofs are
identical, and thus the proof of Theorem \ref{thm:PT-log-concave}
is omitted.

\begin{Theorem} \label{thm:PT-log-concave}
For every $0<\delta<1$ there exist constant $c(\delta)$ for which
the following holds. Let $\mu$ be an isotropic, log-concave measure
on $\R^n$ and let $K \subset \R^n$ be a convex symmetric body. If
$X_1,...,X_k$ are independent, distributed according to $\mu$, then
with probability at least $1-\delta$,
$$
{\rm diam}({\rm ker} \Gamma \cap K) \leq q_k^*(K),
$$
where
$$
q_k^*(K) = \inf\left\{\rho>0 : \rho \geq c(\delta)\gamma_2(K\cap
\rho S^{n-1} ,\psi_2)\sqrt{\frac{\gamma_2(K\cap \rho S^{n-1}
,\psi_2)}{k}} \right\}.
$$
\end{Theorem}
If $\mu$ is a subgaussian measure then for every $A \subset \R^n$,
$\gamma_2(A,\psi_2) \leq c\ell_*(A)$, and thus $\gamma_2(K\cap \rho
S^{n-1} ,\psi_2) \leq \sqrt{n}$. Therefore, Theorem
\ref{thm:PT-log-concave} recovers Theorem \ref{thm:MPT1} up to a
$\sqrt{\log{n}}$ factor. Of course, the bound on the probability is
considerably weaker. On the other hand, Theorem
\ref{thm:PT-log-concave} holds for a much wider family of measures
because the bound given in Theorem \ref{thm:MPT1} depends on the
equivalence constant between the $\psi_2$ and $\ell_2$ metrics
endowed on $\R^n$.

\subsection{Sampling from an isotropic, log-concave measure}

A question which was originally studied in
\cite{KLS,Bour,Rud,GiaMil,GueRud} is the following: how many
points sampled from an isotropic, convex, symmetric body are
needed to ensure that the random operator $k^{-1}\sum_{i=1}^k
\inr{X_i,\cdot}e_i$ is an almost isometric embedding of $\ell_2^n$
in $\ell_2^k$? In other words, that with probability at least
$1-\delta$, for every $\theta \in S^{n-1}$,
$$
1-\eps \leq \frac{1}{k} \sum_{i=1}^k \inr{X_i,\theta}^2 \leq
1+\eps.
$$
\begin{Theorem} \label{thm:best-prev-est} \cite{Rud,GiaMil,GueRud}
For every $0<\eps,\delta<1$ there is a constant $c(\eps,\delta)$
for which the following holds. Let $X_1,...,X_k$ be independent
random variables, distributed according to the volume measure of a
convex, symmetric body in isotropic position. If $k \geq
c(\eps,\delta)n\log^2n$, then with probability at least
$1-\delta$, for every $\theta \in S^{n-1}$,
$$
1-\eps \leq \frac{1}{k} \sum_{i=1}^k \inr{\theta,X_i}^2 \leq
1+\eps.
$$
\end{Theorem}
The estimate of $k \sim n \log^2n $ was first proved by Rudelson
\cite{Rud}. Previously, Bourgain showed \cite{Bour} how to obtain
this result with a slightly weaker estimate of $k \sim n \log^3n$,
but then Giannopoulos and Milman \cite{GiaMil} demonstrated that
Bourgain's method can actually give the same estimate as
Rudelson's.

The proofs of Bourgain and Rudelson use very different arguments.
Rudelson's proof is based on a noncommutative Khintchine
inequality, due to Lust-Piquard and Pisier \cite{LuPi}, namely, a
bound on Rademacher averages of the form $\E \|\sum_{i=1}^k \eps_i
x_i \otimes x_i\|^p_{\ell_2 \to \ell_2}$ for $p \geq 1$, where
$(x_i)_{i=1}^k \in \ell_2$ and $(\eps_i)_{i=1}^k$ are independent,
symmetric, $\{-1,1\}$-valued random variables. The fact that the
set indexing the empirical process is exactly the Euclidean sphere
is essential in the proof and the argument can not be modified to
handle any other indexing sets - not even other subsets of the
sphere.

Bourgain's proof uses a similar technique to the one we used here,
which relies on the following version of Theorem \ref{thm:sets}. The
formulation we present here is from \cite{GiaMil}.
\begin{Lemma} \label{lemma:Bour}
Let $\delta \in (0,1)$ and let $X_1,...,X_k$ be points in $\R^n$
sampled according to an isotropic log-concave measure. If $k \leq
c\delta \exp(\sqrt{n})$ then with probability at least $1-\delta$,
for every $I \subset \{1,...,k\}$,
$$
\left\|\sum_{i \in I} X_i \right\| \leq c_1(\delta)
\left(\sqrt{\log{k}}\sqrt{|I|}\sqrt{n}+|I| \log{k} \right).
$$
In particular, with probability at least $1-\delta$, for every $t
\geq c(\delta)\log{k}$ and every $x \in S^{n-1}$,
\begin{equation} \label{eq:tail1}
\left|\left\{ i : \inr{x,X_i} \geq t\right\} \right| \leq
c_2(\delta)\frac{n\log{k}}{t^2}.
\end{equation}
\end{Lemma}


Bourgain's method was generalized in \cite{GiaMil} in which the
following theorem was established:

\begin{Theorem} \label{thm:GiaMil}
Let $p>0$ and $0 < \delta < 1$. There exists $n_0(\delta)$ such
that for every $n \geq n_0(\delta)$, every log-concave measure on
$\R^n$, every $k \geq k_0(\delta,p)$ and every $\theta \in
S^{n-1}$,
$$
c_p \leq \frac{1}{k}\sum_{i=1}^k |\inr{X_i,\theta}|^p \leq C_p
$$
where $c_p$ and $C_p$ depend only on $p$ and
\begin{equation*}
k_0(\delta,p) =c(\delta,p)
\begin{cases}
n  & \text{if $0<p<1$}, \\
n\log^p{n} & \text{if $1 \leq p \leq 2$}, \\
\min\{(p-2)^{-1},\log{n}\} (n\log{n})^{p/2} & \text{if $p>2$}.
\end{cases}
\end{equation*}
\end{Theorem}
Note that this bound is isomorphic in nature rather than almost
isometric, though in the case $p=2$ the proof of Theorem
\ref{thm:GiaMil} can be modified to give an almost isometric
estimate.

Recently, Gu\'{e}don and Rudelson \cite{GueRud} were able to bound
$$
\E_\eps \sup_{y \in K} \sum_{i=1}^k \eps_i |\inr{x_i,y}|^p,
$$
for any $x_1,...,x_k \in \R^n$, where $K \subset B_2^n$ is a
convex, symmetric body which has a $q$-power type modulus of
convexity. The method of proof is based on majorizing measures,
and can be used to bound $\E\|\mu_k-\mu\|_{F^p}$ for
$F=\{\inr{x,-}:x \in K\}$ as long as $p \geq q \geq 2$. It turns
out that the dominant factor in the bound is
$$
\left(\E \max_{1 \leq i \leq k}\|X_i\|^2 \cdot \E \max_{1 \leq i
\leq k}\|X_i\|^{p-2}_{K^\circ}\right)^{\frac{1}{p}}.
$$
For $K=B_2^n$ this approach yields the best known estimates for
$\E\|\mu_k-\mu\|_{F^p}$ for $p \geq 2$, and the resulting estimate
on the required size of the sample is $k \sim
c(\eps,\delta,p)n^{p/2}\log{n}$, and in particular, for $p=2$
gives the best known estimate of $k \sim c(\eps,\delta)n\log{n}$.
Let us mention that for $p=2$ this result is not helpful for
``small" subsets of the sphere, and the best bound that one can
establish for such subsets coincides with the one obtained for the
whole sphere.

All the known bounds, including \cite{GueRud} and ours, are based on
the behavior of the random variable $\|X\|$. The best estimates on
$\|X\|$ are due to Paouris \cite{Pao}:

\begin{Theorem} \label{thm:Paouris}
There are absolute constants $c_1$ and $c_2$ for which the
following holds. Let $X$ be distributed according to an isotropic
log-concave measure on $\R^n$. Then, for every $p \leq
c_1\sqrt{n}$, $(\E\|X\|^p)^{1/p} \leq c_2 \sqrt{n}$.
\end{Theorem}

Theorem \ref{thm:Paouris} immediately leads to a removal of a
logarithmic factor in \eqref{eq:tail1}, though not to an improved
level of truncation; thus, the estimate of Theorem
\ref{thm:GiaMil} in the case $p=2$ remains unchanged despite the
improved tail estimate.

The properties of an isotropic log-concave measure which will be
used below are that for suitable absolute constants $C$ and $C_1$,
\begin{description}
\item{1.} linear functionals have a subexponential tail - that is, for every $x \in
\R^n$,
$$
\|\inr{X,\cdot}\|_{\psi_1} \leq C\|x\|,
$$
and
\item{2.} By Theorem \ref{thm:Paouris}, for $n \leq k \leq
\exp(c\sqrt{n})$,
$$
\E\max_{1 \leq i \leq k} \|X_i\| = \E\max_{1 \leq i \leq k}
\sup_{x \in S^{n-1}} |\inr{x,X_i}| \leq c_1\sqrt{n}.
$$
\end{description}
Therefore, in light of Theorem \ref{thm:main-trunk}, it is enough
to consider the truncated measure $\nu$ on $\R^n$ which is
supported on a ball of radius $c_2(\delta)\sqrt{n}$ to bound
$\|\mu_k-\mu\|_{F^p}$. The main ingredient in our method is to
bound $\gamma_2(S^{n-1},\nu)$ and to that end we shall estimate
$$
\ell_E = \E\|\sum_{i=1}^n g_ie_i\|_E,
$$
where $g_1,...,g_n$ are standard, independent Gaussian variables.
The particular norm $\| \ \|_E$ we consider is the one endowed on
$\R^n$ by the $\psi_2(\nu)$ structure, formally defined for every
$t \in \R^n$ by $\|t\|_E = \|\inr{t,Y}\|_{\psi_2}$.

\begin{Lemma} \label{lemma:M-estimate}
There exists an absolute constant $c$ for which the following
holds. Let $\nu$ be a probability measure on $\R^n$ and set $Y$ to
be distributed according to $\nu$. If $Z=\|Y\|$ and $E=(\R^n, \| \
\|_{\psi_2})$, then $\ell_{E} \leq c\|Z\|_{\infty}$.
\end{Lemma}

\proof Fix $\rho$ to be named later and consider the gaussian
vector $G=(g_1,...,g_n)$, where $(g_i)_{i=1}^n$ are independent,
standard gaussian random variables. Let $\|Z\|_\infty =D $ and
since $\|f\|_{\psi_2} \leq \E \exp(f^2)$ then
\begin{align*}
\frac{\ell_E}{\rho} \leq & \E_Y \E_g \exp \left(\frac{\sum_{i=1}^n
g_i \inr{e_i,Y}}{\rho} \right)^2.
\end{align*}
Recall that $\sum_{i=1}^n g_i \inr{Y,e_i}$ is distributed as
$g\|Y\|$, and thus,
\begin{align*}
\frac{\ell_E}{\rho} \leq & \E_Y \E_g \left(1+\sum_{m=1}^\infty
\frac{\left(\sum_{j=1}^n g_i \inr{Y,e_i} \right)^{2m}}{m! \rho^{2m}}
\right)
\\
\leq &  \left(1+\sum_{i=1}^\infty \frac{1}{m!} \E|g|^{2m}\left(
\frac{D}{\rho}\right)^{2m}\right)
\\
= & \E \exp \left(\frac{Dg}{\rho } \right)^2 \leq 2,
\end{align*}
if one selects $\rho=cD$. Therefore, $\ell_E \leq cD$, as claimed.
\endproof

\begin{Definition}
For two sets $A,B \subset \R^n$, let $N(A,B)$ be the minimal
number of translates of $B$ needed to cover $A$, that is, the
minimal cardinality of a set $\{x_1,...,x_m\}$ such that $A
\subset \bigcup_{i=1}^m (B+x_i)$.
\end{Definition}
Note that if $B$ is a unit ball of a norm on $\R^n$ then $N(A,\eps
B)$ are the covering numbers of $A$ with respect to the metric
endowed by $B$.


\begin{Corollary} \label{cor:entropy-est}
There exists an absolute constant $c$ such that for every $\eps \geq
1/2$,
$$
\log N(B_2^n, \eps B_E) \leq c\frac{n}{\eps^2},
$$
and for $0<\eps<1/2$,
$$
\log N(B_2^n, \eps B_E) \leq cn\log \left(\frac{1}{\eps}\right),
$$
where $B_E$ is the unit ball of $(\R^n,\psi_2(\nu))$ and $B_2^n$ is
the Euclidean unit ball.
\end{Corollary}

\proof Recall that if $Z=\|Y\|$ then $\|Z\|_{\infty} \leq c_1
\sqrt{n}$ for a suitable absolute constant. By the dual Sudakov
Theorem \cite{PajTom2}, $\log N(B_2^n, \eps B_E) \leq
c_2\ell_E/\eps^2$, and applying Lemma \ref{lemma:M-estimate},
$\ell_E \leq c_3\sqrt{n}$, from which the first part of the claim
follows. Turning to the second part, by a standard volumetric
estimate (see, e.g. \cite{Pis}), and since $B_E$ is a unit ball of
a norm on $\R^n$, $N(\frac{1}{2}B_E,\eps B_E) \leq (1/2\eps)^n$.
Therefore,
$$
N(B_2^n,\eps B_E) \leq N\left(B_2^n, \frac{1}{2}B_E\right) \cdot
N\left(\frac{1}{2}B_E, \eps B_E \right) \leq \exp(c_4n)
\left(\frac{1}{2\eps}\right)^n.
$$
\endproof
Using Corollary \ref{cor:entropy-est} one can bound
$\gamma_2(F,\psi_2(\nu))$ by applying Theorem \ref{thm:talagrand}
for the space $\ell_2^n$ which is $2$-convex, and with $d$ being
the $\psi_2(\nu)$ metric endowed on $\R^n$.

\begin{Corollary} \label{cor:gamma-est}
Let $\mu$ be an isotropic log-concave measure on $\R^n$ and set
$\nu$ to be its truncation as above. Then for $n \leq k \leq
\exp(c_1\sqrt{n})$,
$$
\gamma_2(S^{n-1}, \psi_2(\nu) ) \leq c_2\sqrt{n\log{n}},
$$
where $c_1$ and $c_2$ are absolute constants.
\end{Corollary}
\proof The proof is immediate from Theorem \ref{thm:talagrand},
the entropy estimate in Corollary \ref{cor:entropy-est}, combined
with the fact that for every $\theta \in S^{n-1}$,
$|\inr{Y,\theta}| \leq \|Y\| \leq c\sqrt{n}$, and thus ${\rm
diam}(S^{n-1},\psi_2(\nu)) \leq c\sqrt{n\log{k}} \leq
c\sqrt{n\log{n}}$.
\endproof

Let us remark that we believe this estimate is suboptimal by a
factor of $\sqrt{\log{n}}$.

Combining Corollary \ref{cor:gamma-est} with Theorem
\ref{thm:main-trunk}, we obtain the following (most likely,
suboptimal) estimate of $\|\mu_k-\mu\|_{F^p}$, which we only state
for $p > 2$. This estimate recovers the best known result for
$p>2$, and was originally established in \cite{GueRud}.

\begin{Theorem} \label{thm:p-sphere}
For every $0<\eps$, $0<\delta<1$ and $p >2$ there exists a constant
$c(\eps,\delta,p)$ for which the following holds. With probability
at least $1-\delta$, if $k \geq k_0$,
$$
\sup_{\theta \in S^{n-1}} \left| \frac{1}{k} \sum_{i=1}^k
|\inr{X_i,\theta} |^p - \E|\inr{X,\theta}|^p \right| < \eps,
$$
provided that $k_0 \geq c(\eps,\delta,p)n^{p/2}\log{n}$ for $p>2$.
\end{Theorem}

\proof Let $n \leq k \leq \exp(c_1\sqrt{n})$. Using the notation
of Theorem \ref{thm:main-trunk}, observe that
$\gamma_2(S^{n-1},\psi_2(\nu)) \leq c_2\sqrt{n\log{n}}$, ${\rm
diam}(S^{n-1},\psi_1) \leq c_2$, $H_k \leq c_2\sqrt{n}$. Also, if
$k \geq c_3n\log{n}$, then for $p > 2$, $\theta$ can be taken as
$\theta \leq c_4 \log \log{n}$, from which the claim is evident.
\endproof
Let us remark that if one could select $k \leq cn\log{n}$ it would
be possible to take $\theta$ at the level of an absolute constant.
This would be the case if the logarithmic term in the estimate on
$\gamma_2(S^{n-1},\psi_2(\nu))$ were to be removed and would lead to
the optimal bound for any $p>1$.

\footnotesize {
}

\vskip4cm

\end{document}